%%%%%%%%%%%%%%%%%%%%%%%%%%%%%%%%%%%%%%%%%%%%%%%%%%%%%%%%%%%%%%%%%%%%%%%%%%%%%%

\documentstyle{amsppt}
\nologo
%\NoPageNumber
%
% ------ Macros ------
%
\font\b=cmr10 scaled \magstep4
\def\bigzerou{\smash{\kern-25pt\lower1.7ex\hbox{\b 0}}}
\hsize=360pt
\vsize=540pt
\hbadness=5000
\tolerance=1000
\NoRunningHeads

\def\Q{\Bbb Q}

\def\al{\alpha}

\def\la{\lambda}

\def\Fp{\Bbb F_p}
\def\Fl{\Bbb F_l}
\def\F5{\Bbb F_5}

\def\Z{\Bbb Z}

\def\C{\Bbb C}

\def\CL{\Cal L}

\def\CO{\Cal O}

\def\CM{\Cal M}

\def\no{\noindent}
\def\la{\lambda}

\def\ti{\times}
\def\si{\sigma}

\def\Ga{\Gamma}

\def\ot{\otimes}
\def\part{\partial}

\magnification=\magstep1
%\pageno=27

\topmatter
\title zero-cycles on self-product of modular curves\endtitle
\author
Kenichiro Kimura
\endauthor
\endtopmatter
\vskip -3ex
\noindent\hskip 2em

\vskip 3ex
\CenteredTagsOnSplits
% ------  Document  ------
%
\document
\subhead \S 1 Introduction\endsubhead
 Let $X$ be a projective smooth variety over
$\Q$. Let $S$ be a finite set of primes such that $X$ has 
a projective smooth model $\Cal X$ over $U:=\text{spec}\Z[\frac1S]$.
Then from the works of Bloch(\cite{Bl1}) and Sherman, 
we have the following exact localization sequence in algebraic 
$K$-theory, using Gersten's conjecture for $\Cal K_2$
in a mixed characteristic setting that was proven by Bloch(\cite{Bl2})

$$H^1(\Cal X,\,\Cal K_2)\to H^1(X,\,\Cal K_2)
\overset{\part}\to{\to} \underset{p\in U}\to\oplus Pic(X_p)
\to CH^2(\Cal X)\to CH^2(X)\to 0$$ where $X_p$ is the fiber of
$\Cal X$ at the prime $p$.

A conjecture of Beilinson on the special values of $L$-function
and Tate conjecture tell us that the cokernel of $\part$
is torsion. See the remark after Theorem 2.5 in \cite{Lang}
for an account of this.

In this paper we present a modest evidence for this conjecture.
We consider the case where $X$ is a self-product of a curve $C$ 
which is modular in the sense defined below. 
We make use of elements in the $K$-cohomology
group $H^1(X, \Cal K_2)$ which are constructed by Flach (\cite{Fl})
and Mildenhall(\cite{Mi}). Our main result (theorem 2.3) is that 
the boundary map
$$\part:\,H^1(C\ti C, \Cal K_2)\to 
\bigoplus_{p\in \tilde{S}} Pic(C\ti C\,\text{mod}\,p)$$
has a torsion cokernel. Here $\tilde{S}$ is a set of primes 
which satisfy
certain conditions explained below. By using an idea 
of Murty and Baba in \cite{BaMu} and results by Ribet(\cite{Ri})
and Momose(\cite{Mo}) 
on the image of the Galois representaion associated to $C$,
it is shown in the proposition 3.1 
that under certain mild assumption on $C$ the set $\tilde{S}$ has density 1.

\no When $C$ is an elliptic curve, $\tilde{S}$ contains all but finite
number of primes. 
By using this fact,
Langer and Saito give a finiteness result for the torsion subgroup
of the Chow group $CH^2(C\ti C)$(\cite{LS}). When genus($C)> 1$
this is still not known. So the finiteness result we obtain is 
about $CH^2(C\ti C\ti\,\text{Spec}\, \Q_p)$ for the primes $p\in \tilde{S}$. 
The conditions on the primes are almost equivalent to 
irreducibility of the characteristic polynomial
of the Hecke operator $T_p$ on $\Gamma(C,\Omega^1_C)$. In the appendix
we discuss the relation of this condition with Maeda's conjecture
(\cite{HiMa},\cite{BaMu})
on irreducibility of Hecke polynomials.

Th author thanks Shuji Saito for discussions and Andreas Langer 
for careful reading of earlier version of this paper and
correcting some errors.

\subhead \S 2 The result\endsubhead
For an abelian group
$M$, $M_\Q$ denotes $M\ot \Q$.

\no Let $f$ be an element of $S_2(\Ga_0(N))$ or $S_2(\Ga_1(N))$ for some 
integer $N$ which is an eigenfunction of the Hecke operators 
$T_n$ for all $n>0$, and let $f|T_n=a_nf$. Let $K$ be the subfield
of $\C$ generated over $\Q$ by $a_n$ for all $n$.
Then by the Theorem 7.14 in \cite{Sh} there exists a unique abelian
subvariety $A$ of $J_0(N)$ or $J_1(N)$ depending on 
whether $f\in S_2(\Ga_0(N))$ or $f\in S_2(\Ga_1(N))$
which is defined over $\Q$ with a homomorphism $\theta:\,K\to
End_{\Q}(A)\ot\Q$ which has the following properties:

(1) $\text{dim}\,A=[K:\Q].$

(2) $\theta(a_n)$ is the restriction of $T(n)$ to $A$.

Let $C$ be a projective smooth curve over $\Q$.
\proclaim{Definition} The curve $C$ is modular if 
the jacobian $Jac(C)$ of $C$ is isogenous over $\Q$ to
an abelian variety $A$ which is associated with a modular 
form $f$ in $S_2(\Ga_0(N))$ or $S_2(\Ga_1(N))$ for some $N$ 
in the above sense.\endproclaim For example 
the curves $X_0(23)$ and $X_1(13)$ are modular in this sense.

\no Assume that $C$ is modular. Then there is an algebraic number field
$K$ with $[K:\Q]=g=\text{genus}(C)$ and a homomorphism 
$\theta:\,K\to End_{\Q}(Jac(C))\ot\Q$.
We identify $K$ with $\theta(K)$. 
In the following we denote by $X$ the modular curve $X_0(N)$ or
$X_1(N)$ depending on whether $C$ is modular for $\Ga_0(N)$ or
$\Ga_1(N)$. If $C$ is modular for $\Ga_1(N)$ we assume that 
the modular form $f$ associated to $C$ satisfies the condition
(7.5.9) in \cite{Sh} i.e. the field $K$ is generated by $a_n$
for all $n$ prime to $N$.

\no There is a finite set of primes $S$ which contains all the primes
dividing $N$ such that there are projective smooth models $\Cal X$ 
and $\Cal C$ over $\text{Spec}\, \Z[\frac1S]$ of $X$ and $C$ respectively.
There is also an algebraic correspondence
$\Ga\subset \Cal X\ti \Cal C$ which is finite over $\Cal C$ such that
the generic fiber $\Ga_\Q$ of $\Ga$ induces
a surjective morphism ${\Ga_\Q}_*:Jac(X)\to Jac(C)$. 

\no For a scheme $\Cal Y$ over $\text{Spec}\,\Z[\frac1S]$ and for a prime 
$p\not\in S$ $Y_p$ denotes $\Cal Y\ti_{\text{Spec}\, \Z[\frac1S]} 
\text{Spec}\,\Fp.$

\no  Let $p\not\in S$ be a prime.
First we consider the case where $X=X_0(N).$
Let  $\Cal X_0(Np)$ be the modular curve over
$\text{Spec}\, \Z[\frac1S]$ classifying 
pairs $(E,\,G)$ of elliptic curves and cyclic subgroups of $E$ of order
$Np$.The Hecke correspondence $T(p)\subset \Cal X\ti \Cal X$ 
is the image of $\Cal X_0(Np)$
under the map $(E,\,G)\mapsto (E,\,pG)\ti (E/G',\,G/G')$. Here
 $G'$ denotes the subgroup of $G$ of order $p$.
By (\cite{Sh}, Corollary
7.10) The divisor $T(p)_p$ on $X_p\ti X_p$ 
is the sum of the graph of Frobenius endomorphism and its 
transpose. Since the composition of correspondences
$\Ga\circ {}^t\Ga=(\text{deg}\,\Ga)\Delta\subset \Cal C\ti \Cal C$
where $\text{deg}\,\Ga$
is the degree of $\Ga$ over $\Cal C$,
it follows that the correspondence
 $$(\Ga\circ T(p)\circ {}^t\Ga)_p=\text{deg}\,\,\Ga(\Ga_F+{}^t\Ga_F)$$ 
 where $\Ga_F$
is the graph of frobenius endomorphism $F$ on $C_p$. 

\no Thus we see that there is the equality
$$(\text{deg}\,\Ga)(F^2-a_p F+p)=0 \tag 1$$ 
in $End_{\Fp}(Jac(C_p))_\Q$.

\no Next we consider the case where $X=X_1(N).$
Let $\Cal X_1(N;p)$ be the modular curve over $\text{Spec}\, \Z[\frac1S]$
classifying triples
$(E,P,G)$ of elliptic curves, points of exact order $N$ and cyclic
subgroups of order $p$. The Hecke correspondence 
$T(p)$ on $\Cal X\ti \Cal X$
is the image of $\Cal X_1(N;p)$ by the map $(E,P,G)\mapsto 
(E,P)\ti (E/G,P)$. Then again by (\cite{Sh}, Corollary
7.10) the divisor $T(p)_p$ on 
$X_p\ti X_p$ is equal to $\Ga_{Fr}+\Ga_{<p>}\circ {}^t\Ga_{Fr}$
where $<p>:\,\Cal X\to \Cal X$ is the map 
$(E,P)\mapsto (E,pP)$. We have the equality 
 $$(\Ga_p\circ T(p)\circ {}^t\Ga)_p
 =(\text{deg}\,\,\Ga)\Ga_F+(\Ga\circ \Ga_{<p>}\circ
 {}^t\Ga)_p\circ {}^t\Ga_F.$$
We denote by $<p>'$ the endomorphism of $Jac(C_p)$ induced by
the correspondence $(\Ga\circ \Ga_{<p>}\circ {}^t\Ga)_p$.
Since the morphism $<p>$ is rational
over $\Q$ we see that there is the equality
$$(\text{deg}\,\Ga)(F^2-a_p F)+p<p>'=0 \tag 2$$
 in $End_{\Fp}(Jac(C_p))_\Q$.

We know by Flach, Mildenhall(\cite{Fl},\cite{Mi}) in the case
$X=X_0(N)$ and Weston(\cite{We}) in the case $X=X_1(N)$
that there is an element $\al_p\in H^1(X\ti X, \Cal K_2)_\Q$
such that $\part(\al_p)=\Ga_{Fr}$ the graph of Frobenius endomorphism.
Here the map  
$$ \part:\,\, H^1(X\ti X, \Cal K_2)\to 
\bigoplus_{p\not\in S}Pic(X_p\ti_{\Fp} X_p) \tag 3$$ is the boundary 
map arising from
the localization sequence in algebraic $K$-theory.

\no For a smooth projective variety $V$ over a field $k$, 
The group $H^{n-1}(V, \Cal K_n)$ is isomorphic to the
higher Chow group $CH^n(V,1)$(\cite{Land},\cite{M\"u}).

We define the composition $\Ga\circ \al_p\in CH^2(X\ti C, 1)$ to be
${pr_{13}}_*(pr_{12}^*\al_p \cdot pr_{23}^*\Ga)$
where $pr_{ij}$ are the projections on $X\ti X\ti C$.

 \proclaim{Lemma 2.1} We have the equality 
 $$\part(\Ga\circ \al_p\circ {}^t\Ga)=(\text{deg}\,\Ga)\Ga_F$$
 \endproclaim
 
\demo{Proof} We have the following commutative diagrams:
$$\CD
CH^2(X\ti X\ti C,1)_\Q @>\part>> 
CH^1(X_p\ti X_p\ti C_p)_\Q\\
@Vpr_{23}^*\Ga_\Q\cdot VV        @V pr_{23}^*\Ga_p \cdot VV\\
CH^3(X\ti X\ti C,1)_\Q @>\part >> 
CH^2(X_p\ti X_p\ti C_p)_\Q
\endCD $$

$$\CD
CH^3(X\ti X\ti C,1)_\Q @>\part>> CH^2(X_p\ti X_p\ti C_p)_\Q\\
@V{pr_{13}}_*VV          @V{pr_{13}}_* VV\\
CH^2(X\ti C,1)_\Q @>\part >> CH^1(X_p\ti C_p)_\Q
\endCD $$
Commutativity
of the first diagram follows from the compatibility of the boundary map $\part$
with $K_0(\Cal X^2\ti \Cal C)$ module structure on algebraic $K$-theory. Here 
$\Cal C$ resp. $\Cal X$ is a projective smooth model of $C$ resp. 
$X$ over $\text{Spec}\,\Bbb Z_{(p)}$. See (\cite{Qu}, \S 3) for the definition
of $K_0(\Cal X^2\ti \Cal C)$ module structure.
Commutativity
of the second diagram follows from the functoriality of the map $\part$ on
algebraic $K$-theory with 
respect to the push-forward by projective morphisms. See (\cite{Sri}, (5.11))
for the definition of the push-forward.

\no Hence there is the equality 
$$\align 
& \Ga_p\circ \Ga_{Fr}\\
=&{pr_{13}}_*(pr_{23}^*\Ga_p \cdot \part(pr_{12}^*\al_p))\\
=&{pr_{13}}_*(\part(pr_{23}^*\Ga_\Q \cdot pr_{12}^*\al_p))\\
=&\part({pr_{13}}_*(pr_{23}^*\Ga_\Q \cdot pr_{12}^*\al_p)).
\endalign $$ The same argument for composition with ${}^t\Ga$ concludes
the proof.
\enddemo
By the same argument we also obtain

\proclaim{Lemma 2.2} Assume that $\text{deg}\,\Ga\neq 0.$
Then for an element
 $a\in K\subset End(Jac(C))_\Q$ the composition
$a\circ \Ga_F\in Pic(C_p\ti C_p)_\Q$is contained in $Im(\part)_\Q$.
Here we identify $End(Jac(C))_\Q\subset End(Jac(C_p))_\Q$ 
with a subspace of $Pic(C_p\ti C_p)_\Q.$
\endproclaim

\no From the lemma we see that if $\text{deg}\,\Ga\neq 0$
then $Im(\part)_\Q$ contains
the subring of $End(Jac(C_p))\ot\Q$ which is generated by $F$ and $K$.
 
Since the modular form $f$ associated with $C$ is a Hecke 
eigenfunction for all $T_n$, by the Proposition 3.53 in 
\cite{Sh} $f$ is in $S_2(\Ga_0(N),\,\psi)$ for a unique character 
$\psi:\,(\Z/N\Z )^*\to K^*$. 
Note that $\psi$ is trivial if $f\in S_2(\Ga_0(N)).$

\no As we assume that the condition (7.5.9) in \cite{Sh} holds,
By the theorems 7.14 and 7.16 in loc.cit. it follows
that the characteristic polynomial of the action
of Frobenius endomorphism on the Tate module $T_l(Jac(C))$ is equal to 
$$\prod_{\si:K\hookrightarrow \Bbb C}
(T^2-a_p^\si T+\psi(p)^\si p). \tag 4$$
If this polynomial is separable the dimension the commutator 
of $F$ 

\no in $End_{\Q_l}(T_l(Jac(C)))_\Q$ is equal to 2genus($C$).

\proclaim{Theorem 2.3} Assume that $\text{deg}\,\Ga \neq 0$.
Then the cokernel of the map 
$$ \part:\,\, H^1(C\ti C, \Cal K_2)\to 
\bigoplus_{p}Pic(C_p\ti_{\Fp} C_p) $$
is torsion with $p$ running over the set of primes which satisfy
the following condition:
(1) $p\not\in S$. (2) $F\not\in K$. (3) The polynomial (4) is separable.
\endproclaim
\demo{Proof}
By the Theorem 12.5 in \cite{Mil} the map
$$End_{\Fp}(Jac(C_p))\ot \Q_l\to End_{\Q_l}(T_l(Jac(C_p))_\Q)$$
is injective. So the dimension over $\Q_l$
of the subring of $End_{\Q_l}(T_l(Jac(C_p)_\Q))$
generated by $F$ and $K$ over $\Q_l$ is 2$g$. It follows that 
the subspace $End_{\Fp} (Jac(C_p))_\Q \subset Pic(C_p\ti_{\Fp} C_p)\ot \Q$
is contained in $Im(\part)_\Q$.
The assertion of the theorem follows from the equality
$$Pic(C_p\ti_{\Fp}C_p)\ot \Q=pr_1^*Pic(C_p)_\Q
\oplus End_{\Fp}(Jac(C_p))_\Q \oplus pr_2^*Pic(C_p)_\Q.$$
\enddemo

\remark{Remark} When $C$ is modular for $\Ga_0(N)$ and the field
$K$ is totally real the conditions (2) and (3)
are satisfied if $a_p^\tau\neq a_p^\si$ for any two distinct embeddings
$\si$ and $\tau$ of $K$ into $\Bbb R$ and $\sqrt{p}\not\in K$.\endremark

\proclaim{Corollary 2.4} If a prime $p$ satisfies the conditions
(1), (2) and (3), the non-$p$-primary
torsion part of $CH^2((C\ti C)\ti_\Q \Q_p)$
is finite.
\endproclaim

\demo{Proof} From the works of Bloch(\cite{Bl1}) and Sherman, 
we have the following exact localization sequence in algebraic 
$K$-theory, using Gersten's conjecture for $K_2$
in a mixed characteristic setting that was proven by Bloch(\cite{Bl2})

$$H^1(\Cal Y,\,\Cal K_2)\to H^1(Y,\,\Cal K_2)=CH^2(Y,1)
\overset{\part}\to{\to}  Pic(Y_p)
\to CH^2(\Cal Y)\to CH^2(Y)\to 0.$$
Here $Y=(C\ti C)\ti_\Q \Q_p$, $\Cal Y$ is a projective smooth 
model of $Y$ over  $\text{Spec}\,\Z_p$ and $Y_p=\Cal Y$ mod $p$. By the argument
of Raskind in \cite{Ra, Theorem 1.9}, the reduction map 
$$CH^2(\Cal Y)\{l\}\to CH^2(Y_p)\{l\}$$ is injective for all primes $l\neq p$.
Here $\{l\}$ means the $l$-primary torsion part.
Since $CH^2(Y_p)_{tors}$ is finite(\cite{CTSS}, Theorem 1) the corollary
follows. \enddemo

 \subhead \S 3 About the condition (3) \endsubhead
In this section we consider the condition (3)
in the theorem 2.3. We assume that the curve $C$ is modular
for $\Ga_0(N)$ for some integer $N$. 
In this case the field $K$ is totally real.
Suppose that $J(C)$ satisfies the following condition:
$End(J(C)\ot\bar{\Q})\ot\Q=K$. 

\proclaim{Proposition 3.1} Under the above assumption, we have the inequality
$$ \#\{p\leqq x| a_p\,\,\text{is in a proper subfield of}\,\, K\}
\leqq c\frac{x}{(\log x)^{1+\delta}}$$
for some constants $c$ and $\delta>0$.\endproclaim

\demo{Proof} The proof follows the same argument as in 
\cite{BaMu}. We sketch the outline. Let $L\subset K$ be a proper
subfield. Let $\la\subset \CO_K$ be a prime of degree 
$d\geqq 2$ over $\CL=\la\cap \CO_L$. Let 

$$\rho_{f,\la}:Gal(\bar{\Q}/\Q)\to GL_2(\CO_K/\la)$$
be the Galois representaion $T_l(J(C))\ot_{\CO_K}\CO_K/\la$
where the prime $l$ is the characteristic of $\CO_K/\la$.
The representaion $\rho_{f,\la}$ has the property that
for any prime $p\nmid lN$, we have the equality
$\text{trace}(\rho_{f,\la}(Frob_p))=a_p.$ 
It follows from the proof of the theorem
4.2 in \cite{Mo} that the image of $Gal(\bar{\Q}/\Q)$ in 
$GL_2(\CO_K/\la)=H_\la$ for almost all rational primes $l$
where $H_\la:=\{u\in GL_2(\CO_K/\la)| \det u\in \Fl^*\}.$ 
Note that under our assumtion the field $F$ in loc.cit.is $K$,
the field $L_f$ is $\Q$ 
and the algebra $E=E_f=End_K(H_1(J(C)\ti \C,\,\Q))$.

\no Let $\Cal S\subset  H_\la$ be the set
of elements whose trace lies in the subfield $\CO_L/\CL
\subset \CO_K/\la.$ Then by a counting argument we see that

$$\#\Cal S \leqq e_1N_{\CL}^{2d+1+1/d_2}\quad \text{and}\quad
 \#H_\la =N_\CL^{3d+1/d_2}$$ 
for a positive constant $e_1$. Here $d_2=[\CO_L/\CL:\Fl].$
 
Let $M$ be the fixed field of the kernel of the representation
$\rho_{f,\la}$. By an effective version of the Chebotarev
density theorem in $M$, there exist absolute constants $c_1$ 
and $c_2$ such that if $x\geqq 3$
and $\log x\geqq c_1(\log d_M)(\log\log d_M)(\log\log\log 6d_M)$,
then there is the inequality
$$\#\{p\leqq x | \text{trace}(\rho_{f,\la}(Frob_p))\in \CO_L/\CL\}
\leqq c_2\frac{\#(\Cal S)}{\#(H_\la)}\frac{x}{\log x}.$$
Here $d_M$ is the discriminant of the field $M$. Since we have the 
bound due to Hensel
$$\log d_M\leqq (n_M-1)\sum_{p\in P_L}\log p+n_M\log n_M$$
where $n_M=[M:\Q]$ and $P_L$ is the set of rational
primes which ramify in $L$, we can choose such $\la$ that satisfies
the inequality
$$l\geqq c_3(\log x)^\delta$$ for positive constants $c_3$ and $\delta$.
Hence we obtain the bound
$$\#\{p\leqq x |\text{trace}(\rho_{f,\la}(Frob_p))\in \CO_L/\CL\}
\leqq c_4 \frac{x}{(\log x)^{1+\delta}}$$ for a constant $c_4$.
 Since there are 
only finitely many proper subfields of $K$, this finishes the proof.
\qed \enddemo

\proclaim{Corollary 3.2} Under the same assumption 
as in the proposition  the set of primes which satisfy the conditions
(1), (2) and (3) in the theorem 2.3 has density 1. \endproclaim

\subhead Appendix \endsubhead
We continue to suppose that $C$ is modular for 
$\Ga_0(N).$ The following is a table of 
the genus of the curves $X_0(N)$ and the 
numbers of the primes $p$ such that $p\nmid N,\,p<10000$
and  which do not satisfy
the condition (3) for $X_0(N)$. The curves $X_0(N)$ for
$N=23,41,59$ are modular in our sense.

$$\vbox{ \offinterlineskip
\halign{&\vrule#&\strut\quad#\hfil\quad \cr
\noalign{\hrule}
height2pt & \omit &&\omit &&\omit & \cr
&\hfil $N$ && \hfil genus && \hfil  number of primes & \cr
height2pt & \omit &&\omit && \omit & \cr
\noalign{\hrule}
height2pt & \omit &&\omit && \omit & \cr
& 23 &&
\hfil  2 &&
\hfil 45 & \cr
& 41 &&
\hfil 3    &&
\hfil 1 & \cr
& 59 &&
\hfil 5    &&
\hfil 0 & \cr
height2pt & \omit &&\omit &&\omit & \cr
\noalign{\hrule}
}}$$
The one prime for
$N=41$ is 17.  We used a computer program in 
\cite{St} to compute the characteristic polynomials of the 
Hecke operators. It seems that as the genus of the curve grows,
the primes which do not satisfy the condition (3) become rarer 
very quickly.
We will give a probabilistic argument about reason for
the rareness of the primes which do not satisfy the condition
(3). Let $S'$ be the set of primes 
which do not satisfy the condition (3).
We denote by $R$ the ring $K\cap End(J(C))$. Let $n$ be an integer
which is prime to the discriminant $d_K$ of $K$ and the conductor $F$
of $R$. 
Let $\CM/\text{Spec}\Z[\frac 1{d_KFn}]$ be the moduli space of abelian
schemes of relative dimension genus($C$) with level $n$ structure 
which have real multiplication by $R$. Take a number field
$H$ so that $J(C)_H=J(C)\ti_\Q H$ defines a section
$s_C \in \CM(U)$ for an open set $U$ in Spec$\CO_H$.
For a prime $p\in S'$ the field  $K':=\Q(a_p)$ is a proper subfield
of $K$. Then the characteristic polynomial of $Frob_p$ on $T_l(J(C))$
is equal to 
$$(\prod_{\si:K'\hookrightarrow \Bbb R}
(T^2-a_p^\si T+ p))^n$$
where $n=[K:K']$. This means that there is an isogeny
$$J(C)_p\to (A')^n$$ for an abelian subvariety $A'$ of $J(C)_p$
of dimension $d=[K':\Q]$ which has a real multiplication by the ring
of integers $\CO_{K'}$ of $K'$.
Let $\CL$ be the subscheme of $\CM$ which
corresponds to the abelian schemes which are isogenous to the power
of abelian subschemes with real multiplication by orders in proper
subfields of $K$. The set $S'$ is the image in 
Spec$\Z[\frac 1S]$ of $s_C\cap \CL$. The dimension of $\CM=\text{genus(C)}+1$
and the codimentions of (infinitely many) irreducible components
of $\CL$ are at least $\text{genus(C)}/2$. So as the genus of $C$
grows,
the probability that $s_C$ intersects $\CL$ becomes smaller.

About irredcibility of Hecke polynomials, there is a conjecture 
called Maeda's conjecture which states that the characteristic 
polynomial of the $n$-th Hecke operator $T_n$ on $S_k(SL_2(\Z))$
is irreducible with Galois group $S_d$ for every $n$ (\cite{BaMu},
\cite{HiMa}). We raise a question which is somewhat similar to 
this conjecture:

\proclaim{Question} Is it true that when the genus of the curve $C$
is large, then the characteristic 
polynomial of the $p$-th Hecke operator $T_p$ on $\Gamma(C,\Omega^1_C)$
is irreducible for all but finitely many primes $p$?
\endproclaim
If this question has an affirmative answer, then the condition (3) is satisfied
by almost all primes, and it follows that for all but finitely many
rational primes $p$ the $p$-primary torsion
part of $CH^2(C\ti C)$ is cofinitely generated(cf. \cite{La2}, p278).

Kenichiro Kimura

Institute of Mathematics

University of Tsukuba

Tsukuba, Ibaraki

305-8571

Japan

email: kimurak\@math.tsukuba.ac.jp

\enddocument